\documentclass[11pt]{article} 

\usepackage{pdfsync}
\usepackage[centertags]{amsmath}
\usepackage{amsfonts}
\usepackage{amssymb}
\usepackage{amsthm}
\usepackage{newlfont}
\usepackage[utf8]{inputenc}
\usepackage{mathrsfs}
\usepackage{setspace}
\usepackage{fancyhdr}
\usepackage{epsfig}
\usepackage{graphicx}
\usepackage{accents}

\usepackage{makeidx}
\makeindex

\usepackage[footnotesize,bf]{caption}

\usepackage{verbatim}

\DeclareMathSymbol{\widehatsym}{\mathord}{largesymbols}{"62}

\newtheorem{defn}{Definition}

\newtheorem{prop}[defn]{Proposition}

\newcommand{\pf}{\noindent{\bf Proof }\mbox{   }}

\newcommand{\EV}{\mathbb{E}}    

\newcommand{\ie}{\emph{i.e.}, }
\newcommand{\eg}{\emph{e.g.}, }

\def\midformat{
\setlength{\itemsep}{0pt} \setlength{\parindent}{0mm}
\setlength{\parskip}{0.12in} \setlength{\textheight}{180mm} 
\setlength{\textwidth}{150mm} \setlength{\evensidemargin}{0in}
\setlength{\oddsidemargin}{0in} \setlength{\topmargin}{0in}
\setlength{\hoffset}{1.0cm} \setlength{\voffset}{0.0cm}
\setlength{\headheight}{15pt} \setlength{\headsep}{.5in}
\setlength{\headwidth}{150mm} } \midformat

\newtheoremstyle{theorem}
 {}
 {}
 {\itshape}
 {}
 {\ttfamily}
 {.}
 {.5em}
 {}

\newtheoremstyle{plaintext}
 {}
 {}
 {\upshape}
 {}
 {\ttfamily}
 {.}
 {.5em}
 {}

\theoremstyle{plaintext}

\begin{document}
\makeatother

\pagenumbering{roman} \thispagestyle{empty}

\title{On the Escape of a Random Walk\\From Two Pieces of a Tripartite Set}
\author{Michael Carlisle\footnote{michael.carlisle@baruch.cuny.edu}\\
Baruch College, CUNY}
\maketitle

\begin{abstract}
\noindent Let $\{A, B, C\}$ be a partition of a sample space $\Omega$. For a random walk $S_n = x + \sum_{j=1}^n X_j$ starting at $x \in A$, we find estimates for the Green's function $G_{A \cup B}(x,y)$ and the hitting time $\EV^x(T_C)$ for $x, y \in A \cup B$, with interest in the case where $C$ ``separates'' $A$ and $B$ in a sense (\eg the probability of jumping from $A$ to $B$, or vice versa, before hitting $C$, is small). 
\end{abstract}


\oddsidemargin 0.0in \textwidth 6.0in \textheight 8.5in

\pagenumbering{arabic}
\pagestyle{fancy}
\renewcommand{\headrulewidth}{0.0pt}
\lhead{} \chead{} \rhead{\thepage}
\fancyhead[RO]{\thepage} \fancyfoot{}

\section{Green's functions}

Let $S_n := x + \sum_{j=1}^n X_j$ be a random walk starting at $x$ on a partitioned sample space $\Omega = A \sqcup B \sqcup C$, \ie for any $x, y \in \Omega$, the one-step transition probability is, with $P^x$ the probability measure of the random walk starting at $x$, 
\[ p_1(x,y) = P^x(S_1 = y). \]
Define the first hitting time of $S_n$ on a set $B$ by 
\begin{equation} \label{eq:HittingTimeDefn} 
T_{B} := \inf\{k \geq 0: \, S_k \in B\}. 
\end{equation}
Spitzer, in \cite{Spitzer}, defines the \emph{truncated Green's function}, \index{Green@$G_A$} for $x, y \in A$ of a random walk from $x$ to $y$ before exiting $A$ as the total expected number of visits to $y$, starting from $x$:
\begin{equation} \label{eq:GreenDefn}
G_A(x,y) := \EV^x\bigg[\sum_{j=0}^{\infty} 1_{\{S_j=y; j<T_{A^c}\}}\bigg] = \sum_{j=0}^{\infty} P^x(S_j = y; j < T_{A^c})
\end{equation}
and 0 if $x$ or $y \not \in A$. An elementary result for any random walk (found, for example, in \cite{Spitzer}, or \cite[Sect. 1.5]{LawInt}) is that, for $x, y \in A \subset D$, there are more possible visits inside $D$ than inside $A$:
\begin{equation} \label{eq:GreenComp}
 G_A(x,y) \leq G_D(x,y).
\end{equation}
Starting at a point $x \in A^c$, the \emph{hitting distribution} of $A$ is defined as
\begin{equation} \label{eq:HittingDistDefn}
H_A(x,y) := P^x(S_{T_A} = y).
\end{equation}
The \emph{last exit decomposition} \index{last exit decomposition} of a hitting distribution is based on the Green's function: for $A$ a proper subset of $\Omega$, $x \in A^c$, and $y \in A$, 
\begin{equation} \label{eq:LastExit}
H_{A}(x,y) = \sum_{z \in A^c} G_{A^c}(x,z) p_1(z,y).
\end{equation}
Simple lower bounds for the Green's function $G_{A \cup B}$, by \eqref{eq:GreenComp}, are obvious; for upper bounds for these cases, we examine excursions between $A$ and $B$ before hitting $C$.
\begin{prop} \label{lem:GreenABUB}
For $a,a' \in A$ and $b,b' \in B$, with $\theta_t$ the usual shift operators, 
\begin{align}
T^*_B := & \inf\{k > T_A: S_k \in B\} = T_A + T_B \circ \theta_{T_A}, \notag\\
T^*_A:= & \inf\{k > T_B: S_k \in A\} = T_B + T_A \circ \theta_{T_B}, \notag
\end{align}
and defining
\begin{align}
\psi_a := & \sum_{b' \in B} H_{B \cup C}(a,b') = P^a(T_B < T_C) \label{eq:psix}\\
\sigma_b := & \sum_{a' \in A} H_{A \cup C}(b,a') = P^b(T_A < T_C) \label{eq:sigmax}\\
\rho_a := & \sum_{b' \in B} H_{B \cup C}(a,b') \sigma_{b'} = P^a(T_B, T^*_A < T_C) \label{eq:rhox}\\
\phi_b := & \sum_{a' \in A} H_{A \cup C}(b,a') \psi_{a'} = P^b(T_A, T^*_B < T_C), \label{eq:phix}
\end{align}
we have the Green's function bounds 
\begin{align} 
G_{A}(a,a') \leq G_{A \cup B}(a,a') & \leq G_{A}(a,a') + \frac{\rho_a}{1-\rho_{a'}} G_{A}(a',a') \label{eq:GreenABUBa} \\
G_{B}(b,b') \leq G_{A \cup B}(b,b') & \leq G_{B}(b,b') + \frac{\phi_b}{1-\phi_{b'}} G_{B}(b',b') \label{eq:GreenABUBb} \\
0 \leq G_{A \cup B}(a,b) & \leq \min \left\{\frac{\sigma_b}{1 - \rho_a} G_{A}(a,a), \frac{\psi_a}{1-\phi_b} G_{B}(b,b) \right\}. \label{eq:GreenABUBab} 
\end{align}
\end{prop}
Note that $\psi_a \geq \rho_a$ for every $a \in A$ and $\sigma_b \geq \phi_b$ for every $b \in B$.

\pf We will prove this for \eqref{eq:GreenABUBa} and \eqref{eq:GreenABUBab} (the proof for \eqref{eq:GreenABUBb} matches  \eqref{eq:GreenABUBa}'s proof). By \eqref{eq:GreenDefn}, for $a, a' \in A$, 
\begin{align}
G_{A \cup B}(a,a') & = \sum_{i=0}^{\infty} P^a(S_i = a', i < T_C) \notag\\
 & = \sum_{i=0}^{\infty} [P^a(S_i = a', i < T_C, i < T_B) + P^a(S_i = a', T_B < i < T_C)] \notag\\
 & = G_A(a,a') + \sum_{i=0}^{\infty} P^a(S_i = a', T_B < i < T_C). \label{eq:GAcupB0}
\end{align}
Since $a' \in A$, once the walk enters $B$ it must return to $A$ before hitting $a'$ again. By splitting and switching sums and applying the strong Markov property at $T_B$, 
\begin{align}
G_{A \cup B}(a,a') & = G_A(a,a') + \sum_{i=0}^{\infty} \sum_{b \in B} P^a(S_{T_B} = b, S_i = a', T_B < i < T_C) \notag\\
 & = G_A(a,a') + \sum_{b \in B} H_{B \cup C}(a,b) G_{A \cup B}(b,a'). \label{eq:GAcupB1}
\end{align}

We now switch from \eqref{eq:GreenABUBa} to \eqref{eq:GreenABUBab}: for $G_{A \cup B}(b,a')$, with $b \in B$ and $a' \in A$, decomposing over $A$, and using the strong Markov property at $T_A$, 
\begin{align}
G_{A \cup B}(b,a') & = \sum_{i=0}^{\infty} P^b(S_i = a', i < T_C) \notag\\
 & = \sum_{i=0}^{\infty} \sum_{a'' \in A} P^b(S_i = a', T_A \leq i < T_C, S_{T_A} = a'') \notag\\
 & = \sum_{a'' \in A} H_{A \cup C}(b,a'') G_{A \cup B}(a'',a'). \label{eq:GAcupB2}
\end{align}
We thus have a recurrence relation between \eqref{eq:GreenABUBa} and \eqref{eq:GreenABUBab}.

By the strong Markov property at $T_{a'}$, we have the upper bound 
\begin{equation} \label{eq:InternalPlanarUpperBound}
G_{A}(a'',a') = P^{a''}(T_{a'} < T_{A^c}) G_{A}(a',a') \leq G_{A}(a',a')
\end{equation} 
which yields, by \eqref{eq:sigmax} (for $A \cup B$ instead of $A$), 
\begin{align}
G_{A \cup B}(b,a')  = \sum_{a'' \in A} H_{A \cup C}(b,a'') G_{A \cup B}(a'',a') \leq \sigma_{b} G_{A \cup B}(a',a'). \label{eq:GABzy1}
\end{align}
Combining \eqref{eq:GAcupB1}, \eqref{eq:GABzy1}, and \eqref{eq:rhox} gives us 
\begin{align} 
G_{A \cup B}(a,a') & = G_A(a,a') + \sum_{b \in B} H_{B \cup C}(a,b) G_{A \cup B}(b,a') \notag\\
 & \leq G_A(a,a') + G_{A \cup B}(a',a') \sum_{b \in B} H_{B \cup C}(a,b) \sigma_b \label{eq:GAcupBxy1}\\
 & = G_A(a,a') + G_{A \cup B}(a',a') \rho_a. \notag
\end{align}
In particular, \eqref{eq:GAcupBxy1} gives us 
\begin{align} 
G_{A \cup B}(a',a') \leq \frac{G_A(a',a')}{1 - \rho_{a'}}. \label{eq:GAcupByy}
\end{align}
\eqref{eq:GAcupByy} used again in \eqref{eq:GAcupBxy1} yields \eqref{eq:GreenABUBa}. Proving \eqref{eq:GreenABUBb} similarly, \eqref{eq:GreenABUBb} and \eqref{eq:GAcupByy} applied to \eqref{eq:GABzy1} yields \eqref{eq:GreenABUBab}. \qed

\section{Hitting times}

We now find the expected time of hitting the set $C$, starting from $A$, in terms of hitting $B \cup C$. Lower bounds are simple: just tack the other set on for a quicker hitting time. The upper bounds will require a recursive excursion treatment similar to the proof of Proposition \ref{lem:GreenABUB}.
\begin{prop} \label{lem:ExTCUB}
For $a \in A$ and $b \in B$, defining via \eqref{eq:psix} and \eqref{eq:sigmax}, 
\begin{equation} \label{eq:fAfB}
f_A := \sup_{a \in A} E^a(T_{B \cup C}), \,\,\, f_B := \sup_{b \in B} E^b(T_{A \cup C}), \,\,\, \psi := \sup_{a \in A} \psi_a, \,\,\, \sigma := \sup_{b \in B} \sigma_b,
\end{equation}
we have the expected hitting time bounds 
\begin{align}
E^a(T_{B \cup C}) \leq E^a(T_{C}) \leq E^a(T_{B \cup C}) + \psi_a\left[\frac{f_B  + \sigma f_A}{1 - \psi \sigma}\right] \label{eq:ExTCUBa} \\
E^b(T_{A \cup C}) \leq E^b(T_{C}) \leq E^b(T_{A \cup C}) + \sigma_b\left[\frac{f_A + \psi f_B}{1 - \psi \sigma}\right] \label{eq:ExTCUBb}
\end{align}
\end{prop}

\pf We will prove \eqref{eq:ExTCUBa} (the proof of \eqref{eq:ExTCUBb} is the same). First, decompose $T_C$ along the two possibilities for $T_{B \cup C}$. Recall that $T_{B \cup C} = T_C \iff T_C < T_B$. By the strong Markov property at $T_B$, 
\begin{align}
E^a(T_C) & = E^a(T_C 1_{\{T_{B \cup C} = T_C\}}) + E^a(T_C 1_{\{T_{B \cup C} = T_B\}}) \notag\\
 & \leq E^a(T_{B \cup C}) + \sum_{b \in B} H_{B \cup C}(a,b) E^b(T_C). \label{eq:ExTC1}
\end{align}
Likewise, for $b \in B$, 
\begin{align}
E^b(T_C) \leq E^b(T_{A \cup C}) + \sum_{a' \in A} H_{A \cup C}(b,a') E^{a'}(T_C). \label{eq:ExTC2}
\end{align}
By combining \eqref{eq:ExTC1} and \eqref{eq:ExTC2}, recursing on itself, keeping the first couple terms in terms of $a$, and maximizing the rest via \eqref{eq:psix}, \eqref{eq:sigmax}, and \eqref{eq:fAfB}, we get 
\begin{align*}
E^{a}(T_C) & \leq E^{a}(T_{B \cup C}) + \sum_{b \in B} H_{B \cup C}(a,b) \left(E^b(T_{A \cup C}) + \sum_{a' \in A} H_{A \cup C}(b,a') \left[E^{a'}(T_C)\right]\right) \\
 & \leq E^{a}(T_{B \cup C}) + \sum_{b \in B} H_{B \cup C}(a,b) \left(E^b(T_{A \cup C}) + \sum_{a' \in A} H_{A \cup C}(b,a') \left[ f_A + \psi (f_B + \sigma [...] ) \right]\right),
\end{align*}
which is bounded by 
\begin{align*}
E^{a}(T_C)  & \leq E^{a}(T_{B \cup C}) + \psi_{a} \left(f_B + \sigma [f_A + \psi (f_B + \sigma [...] )]\right) \\
 & = E^{a}(T_{B \cup C}) + \psi_{a} (f_B + \sigma f_A)\sum_{i=0}^{\infty} (\psi \sigma)^i = E^a(T_{B \cup C}) + \frac{\psi_a (f_B + \sigma f_A)}{1 - \psi \sigma}. \qed
\end{align*}

\section{Hitting distributions}

If $y \in A \subset D$, then for $x \in D^c \subset A^c$, we have by \eqref{eq:GreenComp} the monotonicity result
\begin{equation} \label{eq:HittingComp}
\begin{array}{lll}
H_{A}(x,y) & = & \sum_{z \in A^c} G_{A^c}(x,z) p_1(z,y)\\
 & \geq & \sum_{z \in D^c} G_{D^c}(x,z) p_1(z,y) = H_{D}(x,y)
\end{array}
\end{equation}
and the subset hitting time relations (assuming a recurrent random walk) 
\begin{align} 
P^x(T_A = T_D) & = \sum_{z \in A} H_D(x,z); \notag\\
P^x(T_A \neq T_D) & = P^x(T_A > T_D) = \sum_{z \in D \setminus A} H_D(x,z). \label{eq:SubsetHittingComp}
\end{align}
\eqref{eq:HittingComp} and \eqref{eq:SubsetHittingComp} hint at a relationship between the hitting distributions of two sets $C$ and $C \cup A$. We find a bound on this relationship. Let $b \in B$ and $c \in C$. By \eqref{eq:HittingComp} with $D = C \cup A$, there is a probability $p(b, c, C, A)$ such that 
\begin{align}
H_C(b,c) & = H_{C \cup A}(b,c) + p(b,c,C,A). \label{eq:CompareH}
\end{align}
To bound $p(b,c,C,A)$, we rewrite using the definition of $H_C(b,c)$ and decompose along the event $\{T_C < T_A\}$ (whose probability is $1-\sigma_b$ in \eqref{eq:sigmax}):
\begin{align*}
H_C(b,c) & = P^b(S_{T_C} = c) & = P^b(S_{T_C} = c, \, T_C < T_A) + P^b(S_{T_C} = c, \, T_A < T_C); \\
H_{C \cup A}(b,c) & = P^b(S_{T_{C \cup A}} = c) & = P^b(S_{T_{C \cup A}} = c, \, T_C < T_A) + P^b(S_{T_{C \cup A}} = c, \, T_A < T_C).
\end{align*}
Note that 
\[P^b(S_{T_C} = c, \, T_C < T_A) = P^b(S_{T_{C \cup A}} = c, \, T_C < T_A)\]
and 
\[S_{T_{C \cup A}} = c \in C \, \implies \, T_C < T_A,\]
so clearly $P^b(S_{T_{C \cup A}} = c, \, T_A< T_C) = 0$ and we get the simple bound 
\begin{equation} \label{eq:separation}
p(b,c,C,A) = P^b(S_{T_C} = c, \, T_A < T_C)  \leq P^b(T_A < T_C) = \sigma_b.
\end{equation}
If $C$ is a set that ``separates'' $A$ and $B$ in some sense (\eg if the probability distribution of the random walk is based on distance, and $C$ separates $A$ and $B$ into components), then $\sigma_b$ being small reflects the small difference between $H_C$ and $H_{C \cup A}$ (in that it is very likely, starting in $B$, to hit $C$ before $A$).

Note also that $p(C,A)$ is \emph{not} symmetric; \eg $p(A,C) = 1 - p(C,A) = 1 - \sigma_b$.

\singlespacing

\end{document}